\documentclass{amsart}
\usepackage{amsmath,amsfonts,amssymb,amsthm,epsfig,color,tikz, float}
\usepackage[margin=1.3in]{geometry}

\newcommand{\bb}{\mathbb}

\newcommand{\C}{\bb C}

\newcommand{\Z}{\bb Z}
\newcommand{\R}{\bb R}
 
\newcommand{\Q}{\bb Q}

\newcommand{\om}{\omega}

\newtheorem*{obs}{Observation}

\newtheorem*{lemma*}{Lemma}
\newtheorem*{question*}{Question}
\newtheorem*{conj}{Conjecture}

\newtheorem*{theorem*}{Theorem}

\numberwithin{Def}{section}
\numberwithin{equation}{section}

\newtheorem{Theorem}{Theorem}
\numberwithin{Theorem}{section}

\begin{document}
\title[Interval Exchange Maps and Rotations]{Ergodic properties of compositions of Interval Exchange Maps and Rotations}
\author{Jayadev S.~Athreya}
\author{Michael Boshernitzan}
\email{jathreya@iilinois.edu}
\email{michael@rice.edu}
\address{Deptartment of Mathematics, University of Illinois Urbana-Champaign, 1409 W. Green Street, Urbana, IL 61801}
\address{Department of Mathematics, Rice University, Houston, TX.}
\thanks{J.S.A. partially supported by NSF grant
   DMS-1069153.} 
\thanks{M.B. partially supported by NSF grant
   DMS-1102298.} 
\begin{abstract}
We study the ergodic properties of compositions of interval exchange transformations and rotations. We show that for any interval exchange transformation $T$, there is a full measure set of $\alpha \in [0, 1)$ so that $T \circ R_{\alpha}$ is uniquely ergodic, where $R_{\alpha}$ is rotation by $\alpha$.
\end{abstract}
\maketitle

\section{Introduction}\label{sec:intro}

An \emph{interval exchange transformation} (IET) is given by cutting the interval $[0, 1]$ into $m$ subintervals of lengths given by a vector $\lambda = (\lambda_1, \ldots, \lambda_m) \in \R_+^m$, $\sum \lambda_i = 1$ (we denote this set of vectors $\Delta^m$) and a permutation $\pi \in S_m$. The transformation $T = T_{\lambda, \pi} : [0, 1] \rightarrow [0,1]$ is given by gluing the subintervals together in the order given by $\pi$ and preserving the orientation. These are natural generalizations of circle rotations (which can be viewed as exchanges of two intervals), and are closely related to flows on flat surfaces and billiards in Euclidean polygons. IETs preserve Lebesgue measure. 

It was conjectured by Keane, and proved by Masur~\cite{Masur1} and Veech~\cite{Veech} independently that for an irreducible permutation (that is, one that does not fix the 
set $\{1, \ldots, k\}$ for any $k <m$) $\pi$, that $T_{\lambda, \pi}$ is in fact uniquely ergodic  (that is, Lebesgue measure is the \emph{only} preserved measure) for almost all 
$\lambda \in \Delta^m$. 

Boshernitzan \cite{Bosh1} exhibited a Diophantine condition
(Property P) for unique ergodicity of IETs  (see also \cite{Veech3}, an improvement by Veech).
Since this condition is generic (i.e., holds for Lebesgue almost all parameters $\lambda$), 
it provides an alternative approach to Keane's conjecture. 
The condition also allows to establish unique ergodicity of IETs in some 
special situations (e.g., when  $T$  is minimal and all parameters $\lambda_m$  lie in
a quadratic number field, see \cite{Bosh3}). The condition can be applied in a more general
setting \cite{Bosh2} of symbolic flows.

A natural follow-up question is to understand which subsets of IET space inherit the property that almost every IET is uniquely ergodic.  In this note, we study the ergodic properties of the the set of IETs given by compositions of a fixed IET $T$ with an arbitrary rotation $R_{\alpha}$. A simplified version of our main result is:

\begin{Theorem}\label{theorem:iet:ralpha}
Let $T$ be an IET. Then for almost every $\alpha \in [0, 1)$, $T \circ R_{\alpha}$ is uniquely ergodic.

\end{Theorem}

\noindent\textbf{Remarks:}
\begin{enumerate} 
\item The above result holds for \emph{every} interval exchange $T$, including non-ergodic ones (such as the identity). Theorem \ref{theorem:iet:ralpha} provides explicit examples of families of IETs which intersect the set of uniquely ergodic IETs in sets of full measure.
\medskip
\item One consequence of unique ergodicity is the fact that ergodic averages converge for every orbit. In fact, using results of Athreya-Forni~\cite{AtF} this statement can be strengthened to give an upper bound for the rate of convergence for ergodic averages, depending only on the choice of $T$, see Theorem~\ref{theorem:iet:dev}.
\end{enumerate}

\subsection{Plan of Paper}\label{subsec:plan}
In \S\ref{subsec:main}, we state a more general theorem Theorem~\ref{theorem:iet:ralpha:general}. We place our result in the context of recent developments in \S\ref{subsec:hist}.  We recall the required background on flat surfaces and quadratic differentials in \S\ref{sec:qd}, in particular focusing on the construction of surfaces associated to IETs. We prove Theorems~\ref{theorem:iet:ralpha} and~\ref{theorem:iet:ralpha:general} in \S\ref{sec:theoremone}. 
\section{Main Theorem}\label{subsec:main}
\noindent Theorem~\ref{theorem:iet:ralpha} is in fact a special case of the following more general result:

\begin{Theorem}\label{theorem:iet:ralpha:general}
Let $T_1, \ldots, T_k$ be IETs, and $c_1, \ldots, c_k\in\R^{+}$ be positive real numbers. Then for almost every $\alpha \in [0, 1)$, the map $$S_{\alpha} = T_k \circ R_{c_k \alpha} \circ T_{k-1} \circ R_{c_{k-1} \alpha}\circ \ldots \circ T_2 \circ R_{c_2 \alpha} \circ T_1 \circ R_{c_1 \alpha}$$ is uniquely ergodic.

\end{Theorem}

\subsection{Deviation of Ergodic Averages} One of the main consequence of unique ergodicity is control of the ergodic averages for \emph{every} point. Suppose $S:[0, 1] \rightarrow [0, 1]$ is uniquely ergodic with respect to Lebesgue measure and that $f$ is a continuous function. Then for all $x_0 \in [0, 1]$, we have $$\sum_{i=0}^{N-1} f(S^i x_0) = N \int_0^1 f(x) dx + o(N).$$ A general question is whether the $o(N)$-error term can be improved. In our setting, we have:

\begin{Theorem}\label{theorem:iet:dev} Fix notation as in Theorem~\ref{theorem:iet:ralpha:general}. Then there is a $0 \le \beta <1$ depending only on the combinatorics of $T_1, \ldots, T_k$ so that for almost every $\alpha \in [0, 1)$ and any smooth function $f$ on $[0,1]$, 
$$\sum_{i=0}^{N-1} f\left( (S_{\alpha})^i x_0\right) = N \int_0^1 f(x) dx + o(N^{\beta}),$$ for any $x_0$ whose forward orbit is well-defined.
\end{Theorem}

\subsection{Hausdorff dimension of non-ergodic maps} Our proof of Theorem~\ref{theorem:iet:ralpha:general} allows us also to leverage estimates of Masur~\cite{Masurhdim} to conclude:
\begin{Theorem}\label{theorem:hdim} Fix notation as in Theorem~\ref{theorem:iet:ralpha:general}. Then, letting Hdim denote Hausdorff dimension, $$\mbox{Hdim}\,\{\alpha\colon S_{\alpha} \mbox{ is not uniquely ergodic}\} \le \frac 1 2.$$

\end{Theorem}

\subsection{History and Prior Results}\label{subsec:hist}

The study of interval exchange maps and their ergodic properties is an extremely active area. For a beautiful introduction to the combinatorics of interval exchange maps and connections to the study of flat surfaces and Teichm\"uller geodesic flow, see~\cite{Yoccoz}. 

The first main results on the ergodic properties were due to Keane~\cite{Keane}, who proved minimality (every orbit is dense) for IETs satisfying what is known as the infinite distinct orbit condition (i.d.o.c). As discussed above, he conjectured that almost every irreducible interval exchange was uniquely ergodic. The restriction to almost every is required, since there are examples of minimal, non-uniquely ergodic IETs, due to Keynes and Newton~\cite{KeynesNewton}, Keane himself~\cite{Keane2}, and Veech~\cite{Veech1}. This phenomenon does not occur for rotations, as irrationality of the rotation angle implies both minimality and unique ergodicity. More recently, Avila and Forni~\cite{AF}  have shown that almost every interval exchange transformation is weak-mixing, which further emphasizes the contrast with the setting of rotations. 

Our results are somewhat orthogonal to this development. Our proofs rely on results on flows on flat surfaces and billiards. Our main tool in the proof of Theorem~\ref{theorem:iet:ralpha:general} will be the beautiful paper of Kerckhoff-Masur-Smillie~\cite{KMS}, in which it is shown that for every holomorphic quadratic differential, and almost every direction $\theta \in [0, 2\pi)$, the rotation of the vertical foliation by $\theta$ is uniquely ergodic. Similarly, to prove Theorem~\ref{theorem:iet:dev}, we will use the main result of Athreya-Forni~\cite{AtF}, which controls the deviation of ergodic averages for the vertical flow in almost every direction.

\section{Flat Surfaces and Quadratic Differentials}\label{sec:qd}
As discussed in \S\ref{subsec:hist}, there is a close connection between IETs and geodesic flows on singular flat surfaces. We recall some basic definitions and notations. An excellent reference for this material is~\cite{MasurTab}.  Let $\Sigma_g$ be a topological surface of genus $g\geq 1$. 
Let $\Omega_g$ denote the space of holomorphic differentials on $\Sigma_g$, that is, the space
of pairs $(X, \om)$, where $X$ is a genus $g$ Riemann surface and $\omega$ is a holomorphic one-form, that is, a tensor of the form $f(z)dz$ in local coordinates. Any holomorphic differential $\omega \in \Omega_g$ (we drop the $X$ for ease of notation) determines a unique flat metric with conical singularities at the zeros of the holomorphic differential. Given $\omega \in \Omega_g$, one obtains (via integration of the form) an atlas of charts to $\C \cong \R^2$, with transition maps of the form $z \mapsto  z + c$. Vice-versa, given such an atlas of charts, one obtains a holomorphic differential by pulling back the form $dz$ on $\C$.

\begin{figure}[h!]
\begin{center}

\includegraphics[width=0.6\textwidth, height = 25mm]{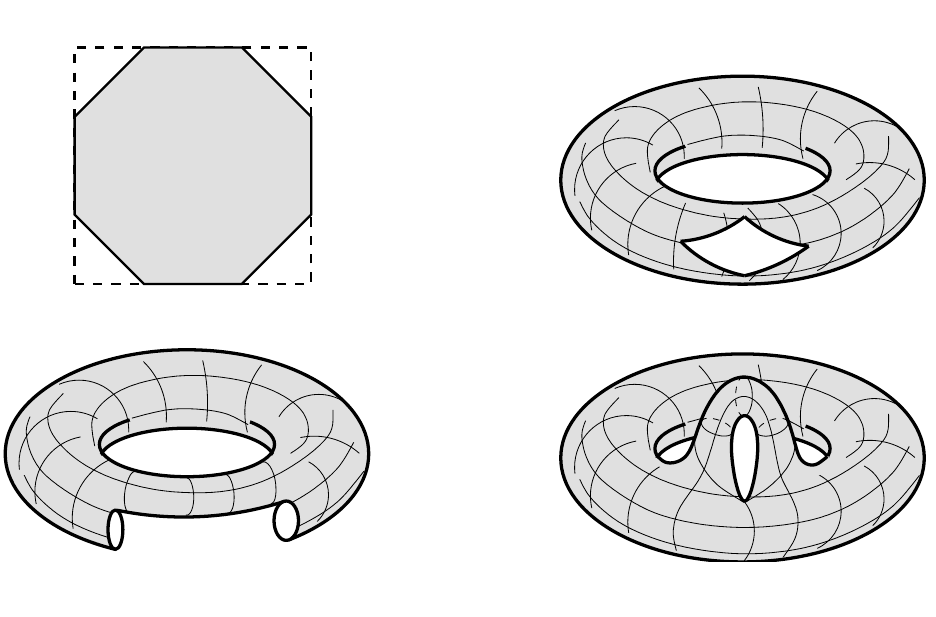}

\caption{Obtaining a flat surface from identifying opposite sides of an octagon. This is a genus $2$ surface with one-form coming from $dz$, which under identifications has one singular point, of order $2$.}
\label{figure:octagon}
\end{center}
\end{figure}

 Any differential $\omega\in \Omega_g$ determines a pair of transverse oriented measured foliations, defined by $\{\operatorname{Re} (\omega) =0\}$ \emph{(vertical foliation)}, $\{\operatorname{Im} (\omega) = 0\}$ \emph{(horizontal foliation)}. These foliations have saddle-like singularities (possibly degenerate) at the zeros of the holomorphic differential, and there are flows associated to moving along leaves of the foliations at unit speed. These flows preserve the Lebesgue measure on the surface arising from the flat metric. Using methods from Teichm\"uller theory, Kerchkoff-Masur-Smillie proved:
 
 \begin{Theorem}\label{theorem:kms} Let $\omega \in \Omega_g$. Then for Lebesgue almost every $\theta \in [0, 2\pi)$, the vertical flow $\{\varphi^t_{\om, \theta}\}_{t \in \R}$ (that is, the flow associated to the vertical foliation) of $e^{i\theta}\omega$ is uniquely ergodic, that is, there is a unique transverse invariant measure.
\end{Theorem}

\noindent The relevance of this theorem to our setting is given by the following observation:

\begin{obs} The first return map of the vertical flow to a transverse interval is an interval exchange transformation. \end{obs}

\noindent Thus, Theorem~\ref{theorem:kms} gives one-parameter families of IETs for which almost every member of the family is uniquely ergodic. We will prove Theorem~\ref{theorem:iet:ralpha} by constructing an appropriate surface and applying Theorem~\ref{theorem:kms}. Similarly, to prove Theorem~\ref{theorem:iet:dev}, we will use the following result, which can be viewed as a quantitative analogue of Theorem~\ref{theorem:kms} (and in fact, where Theorem~\ref{theorem:kms} relies on recurrence estimates for Teichm\"uller geodesic flow, Theorem~\ref{theorem:af} relies on \emph{quantitative recurrence} estimates).

\begin{Theorem}\label{theorem:af}~\cite[Theorem 1.1]{AtF} Let $\omega \in \Omega_g$. Then there is a constant $0 \le \beta \le 1$, depending only on the orders of the zeros of $\omega$, and a function $K_{\omega}: S^1 \rightarrow \R^+$ so that for almost every $\theta \in S^1$, any smooth function $f: X \rightarrow  \R$, and every $x_0$ whose forward orbit is well-defined, $$\left| \int_{0}^T f(\varphi^t_{\om, \theta}(x_0)) dt -T \int_X f dm \right| \le K_{\omega}(\theta) T^{\beta}.$$ Here $m$ denotes the Lebesgue measure associated to the flat metric on $(X, \om)$.
\end{Theorem}

\section{Surface constructions}\label{sec:theoremone}
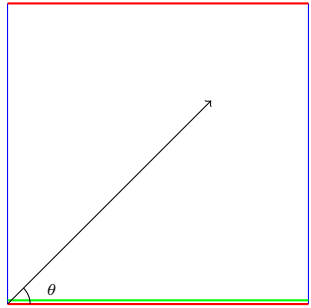
\begin{figure}\caption{The square $[0,1]^2$, with \textcolor{red}{red sides}  identified by the IET $T$, and the \textcolor{blue}{blue sides} are identified by translation. The flow in the direction $\theta$, with inverse slope $\cot \theta =\alpha$, gives $T \circ R_{\alpha}$ as first return map to the horizontal transversal $[0, 1) \times \{0\}$ (marked in \textcolor{green}{green}).\medskip}\label{fig:square}
\vspace{3mm}
\begin{tikzpicture}[scale=1]

\draw[blue] (-2, -2)--(-2, 2);
\draw[blue] (2, -2)--(2, 2);
\draw[thick, red] (-2, -2)--(2, -2);
\draw[thick, red] (-2, 2)--(2, 2);
\draw[thick, green](-2, -1.95)--(2, -1.95);
\draw[->] (-2, -2)--(1/1.414, 1/1.414);
\draw (-1.7, -2) arc(0:45:.3);
\path(-1.6, -1.8)node[right]{\tiny $\theta$};

\end{tikzpicture}
\end{figure}
\medskip

\subsection{Proof of Theorem~\ref{theorem:iet:ralpha}} Fix an IET $T$, and we construct a surface $S_T$ as follows. Consider the square $[0, 1]^2$, with the top identified to the bottom via the IET $T$, and the left and right hand side identified by translation (see Figure~\ref{fig:square}). The surface comes equipped with the differential $dz$, which is preserved by translations.

\noindent The first return map for the flow in direction $\theta$ (with $\cot \theta = \alpha$) to the horizontal transversal is $T \circ R_{\alpha}$. This is the vertical flow for the differential $e^{i\theta} dz$ on the surface $S_T$. Applying Theorem~\ref{theorem:kms}, we obtain Theorem~\ref{theorem:iet:ralpha}.\qed

\begin{figure}\caption{The rectangles $\mathcal R_1, \ldots \mathcal R_k$. The \textcolor{blue}{blue sides} of each $\mathcal R_i$, which have height $c_i$, are identified with each other, and the top of $\mathcal R_i$ is identified with the bottom of $\mathcal R_{i+1}$ via the IET $T_i$. We view the indices cyclically, so $\mathcal R_{k+1} = \mathcal R_k$. Here \textcolor{green}{$T_1$} identifies the \textcolor{green}{green sides}, \textcolor{orange}{$T_{k-1}$} identifies the \textcolor{orange}{orange sides}, and  \textcolor{red}{$T_k$} identifies the \textcolor{red}{red sides}. As above, $\cot \theta = \alpha$.\medskip}\label{fig:rectangles}
\vspace{3mm}
\begin{tikzpicture}[scale=.5]

\draw[blue] (-4, -2)--(-4, 2);
\draw[blue] (0, -2)--(0, 2);
\draw[red] (-4, -2)--(0, -2);
\draw[green] (-4, 2)--(0, 2);
\draw[->] (-4, -2)--(-2+1/1.414, 1/1.414);
\draw (-3.7, -2) arc(0:45:.3);
\path(-3.6, -1.8)node[right]{\tiny $\theta$};
\path(0, 0)node[blue, left]{\tiny $c_1$};
\path(-2, -1)node{$\mathcal R_1$};

\draw[xshift=1cm][blue] (0, -2)--(0, 3);
\draw[xshift=1cm][blue] (4, -2)--(4, 3);
\draw[xshift=1cm]
[green] (0, -2)--(4, -2);
\draw[xshift=1cm]
[yellow] (0, 3)--(4, 3);
\path(5, 0)node[blue, left]{\tiny $c_2$};

\path(3, -1)node{$\mathcal R_2$};

\path(7.5,0)node{\huge $\ldots \ldots$};

\draw[xshift=10cm][blue] (0, -2)--(0, 1);
\draw[xshift=10cm][blue] (4, -2)--(4, 1);
\draw[xshift=10cm]
[purple] (0, -2)--(4, -2);
\draw[xshift=10cm]
[orange] (0, 1)--(4, 1);

\path(12, 0)node{$\mathcal R_{k-1}$};

\path(14, -1)node[blue, left]{\tiny $c_{k-1}$};

\draw[xshift=15cm][blue] (0, -2)--(0, 4);
\draw[xshift=15cm][blue] (4, -2)--(4, 4);
\draw[xshift=15cm]
[orange] (0, -2)--(4, -2);
\draw[xshift=15cm]
[red] (0, 4)--(4, 4);
\path(17, -1)node{$\mathcal R_k$};
\path(19, 0)node[blue, left]{\tiny $c_k$};

\end{tikzpicture}
\end{figure}
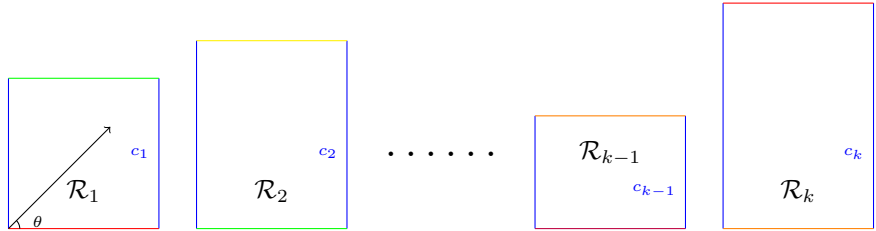

\subsection{Proof of Theorem~\ref{theorem:iet:ralpha:general}} To prove Theorem~\ref{theorem:iet:ralpha:general}, we generalize the above construction. Fix IETs $T_1, \ldots, T_k$, and positive integers $c_1, \ldots, c_k$. Consider $k$ rectangles $\mathcal R_1, \ldots \mathcal R_k$, all of width $1$, and of heights $c_1, \ldots, c_k$. Then glue the top of the first to the bottom of the second with $T_1$, the top of the second to the bottom of the third by $T_2$, and so on in cyclic order, finally gluing the top of the $k^{th}$ to the bottom of the first by $T_k$. Then the flow in direction $\theta = \cot^{-1} \alpha$ will yield as first return map to the horizontal transversal given by the bottom of the first rectangle the map $$S_{\alpha} = T_k \circ R_{c_k \alpha} \circ T_{k-1} \circ R_{c_{k-1} \alpha} \ldots T_2 \circ R_{c_2 \alpha} \circ T_1 \circ R_{c_1 \alpha}.$$ As above, applying Theorem~\ref{theorem:kms}, to directional flows on the surface, we obtain our result.\qed

\subsection{Proof of Theorem~\ref{theorem:iet:dev}} To prove Theorem~\ref{theorem:iet:dev}, we need to pass from the estimate in Theorem~\ref{theorem:af} for deviation of ergodic averages for the flow $\varphi_{\alpha}$ in direction $\alpha$ on the surface $(X, \om)$ constructed above  to an estimate for the deviation of ergodic averages for first return map $S_{\alpha}$. 

A general argument for this procedure was communicated to us by G.~Forni. We can view the flow $\varphi_{\alpha}$ as a suspension flow over the map $S_{\alpha}$ with roof function $\rho_{\alpha}: [0, 1] \rightarrow \R^+$, so that $\int_0^1 \rho_{\alpha}(x) dx = 1$. Note that as $\alpha$ varies, the roof function varies as well. That is, we identify the surface $X$ with the space $$\{ (x, t): x \in [0, 1], 0 \le t < \rho_{\alpha}(x)\}/ \sim,$$ where $(x, \rho_{\alpha}(x)) \sim (S_{\alpha}(x), 1)$. The Lebesgue measure $dm$ on $(X, \om)$ is given by $dx dt$. Given a smooth mean-zero function $f:[0, 1] \rightarrow \R$, and multiply this times a small bump function in the flow direction with unit total
mass to get a smooth mean-zero function $F$ on the surface. 
By Theorem~\ref{theorem:af}, we have, for $T>>0$, $ F(\varphi_{\alpha}(x_0)) dt  = o(T^{\beta})$ and we would like to show, that for $N>>0$,  $\sum_{i=1}^N f(S_{\alpha}^i(x_0)) = o (N^{\beta}).$ We can write (up to a uniformly bounded error) $$\sum_{i=1}^N f(S_{\alpha}^i(x_0)) = \int_{0}^{T_N} F(\varphi_{\alpha}(x_0)) dt,$$ where $T_N =  \sum_{i=0}^{N-1} \rho_{\alpha}(S_{\alpha}^i (x_0)).$ By unique ergodicity of $S_{\alpha}$, we have that, as $N \rightarrow \infty$, for all $x_0$, $$\frac{T_N}{N} \rightarrow \int_0^1 \rho_{\alpha}(x) dx =1.$$ Thus, for $N >> 0$, we have, as desired $$\sum_{i=1}^N f(S_{\alpha}^i(x_0)) = \int_{0}^{T_N} F(\varphi_{\alpha}(x_0)) dt = o(T_N^{\beta}) = o(N^{\beta}).$$  \qed\medskip

\subsection{Proof of Theorem~\ref{theorem:hdim}} Masur~\cite[Main Theorem]{Masurhdim} showed that for any differential $\omega$, $$\mbox{Hdim}\{\theta: e^{i\theta}\omega \mbox{ has non-uniquely ergodic vertical flow}\} \le \frac{1}{2}.$$ Applying this to our surface construction, we have our Theorem~\ref{theorem:hdim}.\qed\medskip

\subsection{Lattice surfaces} We note that the construction of surfaces described above does not give the full richness of possible translation surface structures. In particular, there is an $SL(2, \R)$-action on the space of translation surfaces which gives a sort of \emph{renormalization} dynamics for flat surface flows and interval exchange maps. The stabilizer $SL(X, \om)$ of a surface $(X, \om)$ is known as the \emph{Veech group}, and there are `highly symmetric' surfaces, known as \emph{lattice surfaces}, for which $SL(X, \om)$ is a lattice. A condition for being a lattice surface is that in any direction where the flow is periodic, that the moduli of cylinders in that direction are all commensurable (over $\Q$). For surfaces that arise via our construction, the heights of cylinders in any direction are constant, thus, the condition of moduli of cylinders being commensurable reduces to condition that the lengths of the intervals of the associated IET are rational. This implies that the associated surface must be \emph{square-tiled}, and by results of Gutkin-Judge~\cite{GJ}, the Veech group must be commensurate to $SL(2, \Z)$.

%
%

\section{Further Questions}

\subsection{Negative rotation angles} A natural generalization of our construction would be to consider arbitrary values of $c_i \in \R$. Our construction corresponds to the case where all the $c_i$'s have the same sign. We conjecture:

\begin{conj} Let $T_1, \ldots, T_k$ be IETs, and let $c_1, \ldots, c_k \in \R$ with $\sum c_i \neq 0$. Then for almost every $\alpha \in [0, 1)$, the map $$S_{\alpha} = T_k \circ R_{c_k \alpha} \circ T_{k-1} \circ R_{c_{k-1} \alpha}\circ \ldots \circ T_2 \circ R_{c_2 \alpha} \circ T_1 \circ R_{c_1 \alpha}$$ is uniquely ergodic.

\end{conj}

\noindent It seems natural to attempt to mimic our proof of Theorem~\ref{theorem:iet:ralpha:general}. However, there does not seem to be an obvious construction of a quadratic differential on an orientable surface that allows us to apply Theorem~\ref{theorem:kms}. For all but a countable set of $\alpha$, the transformation $S_{\alpha}$ is minimal, which follows from checking Keane's i.d.o.c.~\cite{Keane}.

\subsection{Weak mixing} It is also natural to consider the question of weak mixing along these families of IETs. In particular, if an IET $T$ is not of rotation type, we conjecture:
\begin{conj} For almost every $\alpha \in [0, 1)$, $T \circ R_{\alpha}$ is weak mixing.
\end{conj}

\subsection{Acknowledgements} J.S.A. would like to thank Rice University for their hospitality for a visit in October 2012 where this work was completed. M.B. would like to thank the University of Illinois for their hospitality for a visit in April 2012 when this work was initiated. We would like to thank Giovanni Forni and Yaroslav Vorobets for useful discussions, and Anton Zorich for his permission to use Figure 1.

\end{document}